\documentclass[12 pt, reqno]{amsart}
\usepackage{ucs}
\usepackage{amsmath}
\usepackage{amsfonts}
\usepackage{amssymb}
\usepackage{hyperref}
\usepackage{url}
\usepackage{amsthm,geometry,latexsym,tabularx,shapepar}
\usepackage{shuffle}
\usepackage[english]{babel}
\usepackage[latin1]{inputenc}
\usepackage{fullpage}
\usepackage[all,2cell]{xy}
\usepackage{graphics,graphicx}
\usepackage[vcentermath]{youngtab}
\usepackage{mathtools}
\usepackage{mathrsfs}
\usepackage{shuffle}
\usepackage{ulem}
\usepackage{cases}
\usepackage{tikz}
\usetikzlibrary{matrix,arrows}
\usepackage{stmaryrd} 
\usepackage{todonotes}
\usepackage{calligra}

\def\C{{\mathbb C}}

\def\N{{\mathbb N}}

\def\Z{{\mathbb Z}}
\def\al{\alpha}

\def\scal#1#2{\langle #1 | #2 \rangle}

\def\LDIAG{\mathbf{LDIAG}}

\gdef\stuffle{\;%
  \setlength{\unitlength}{0.0125cm}%
  \begin{picture}(20,10)(220,580) 
  \thinlines 
  \put(220,592){\line( 0,-1){ 10}} 
  \put(220,582){\line( 1, 0){ 20}} 
  \put(240,582){\line( 0, 1){ 10}} 
  \put(230,592){\line( 0,-1){ 10}} 
  \put(225,587){\line( 1, 0){ 10}} 
  \end{picture}\; 
}

\newtheorem{example}{Example}[section]

\newtheorem{note}{Note}[section]

\newtheorem{theorem}{Theorem}[section]

\newtheorem{lemma}{Lemma}[section]

\def\cqfd{\hfill $\Box$}

\definecolor{MyDarkBlue}{rgb}{0,0.08,0.4}

\def\adots{\mathinner{\mkern2mu\raise1pt\hbox{.}
\mkern3mu\raise4pt\hbox{.}\mkern1mu\raise7pt\hbox{.}}}

\def\up#1{\raise 1ex\hbox{\footnotesize#1}}

\def\mref#1{(\ref{#1})}

\def\adots{\mathinner{\mkern2mu\raise1pt\hbox{.}
\mkern3mu\raise4pt\hbox{.}\mkern1mu\raise7pt\hbox{.}}}

\def\up#1{\raise 1ex\hbox{\footnotesize#1}}

\def\mref#1{(\ref{#1})}

\def\ra{\rightarrow}

\def\path{\rightsquigarrow}

\def\scal#1#2{\langle #1 | #2 \rangle}

\def\LDIAG{\mathbf{LDIAG}}

\newcommand{\eq}[1][r]
   {\ar@<-3pt>@{-}[#1]
    \ar@<-1pt>@{}[#1]|<{}="gauche"
    \ar@<+0pt>@{}[#1]|-{}="milieu"
    \ar@<+1pt>@{}[#1]|>{}="droite"
    \ar@/^2pt/@{-}"gauche";"milieu"
    \ar@/_2pt/@{-}"milieu";"droite"}

\title[Finite Decomposition Semigroups]{Finite Decomposition Semigroups}

\author{Matthieu Deneufch\^atel, G\'erard H. E. Duchamp}
\address{LIPN - UMR 7030 du CNRS, 99, avenue Jean-Baptiste Cl\'ement, Universit\'{e} Paris 13, 93430 Villetaneuse, France}
\email{{matthieu.deneufchatel,ghed}@lipn.univ-paris13.fr}

\date{}

\begin{document}
\begin{abstract}In this paper, we explain the importance of finite decomposition semigroups and present two theorems related to their structure.
\end{abstract}
\maketitle
\tableofcontents

\section{Introduction}
Theories of ``special sums'' have highlighted different products over the indices. For example, Chen's lemma states that the product of two iterated integrals is ruled out by the shuffle product defined by
\begin{equation}
        \begin{aligned}
                1 \shuffle w = w \shuffle 1 & = w ~; \\
                (a u) \shuffle (b v) & = a(u \shuffle (bv)) + b ( (au) \shuffle v) 
        \end{aligned}
\end{equation}
for all words $u, \, v, \, w \in A^*$ and all letters $a,b$ of the alphabet $A$.\\
Indeed (see \cite{JGJY1}), if $\mathscr H$ is a vector space of integrable functions over $(c_1,c_2)$ and $f_1 , \dots , f_n$ some functions of $\mathscr H$, define the following integral:
\begin{equation}
   \langle f_1 \dots f_n \rangle = \int_{c_1}^{c_2} dy_1 \int_{c_1}^{y_1} \dots \int_{c_1}^{y_{n-1}} d y_n \, f_1(y_1) \dots f_n(y_n)
\end{equation}
(considered as a linear form defined on ${\mathscr H}^{\otimes n}$). 

If the functions $\phi_{a_i}$ are indexed by letters of the alphabet $A$, we associate to $w = a_{i_1} \dots a_{i_{|w|}}$ the integral
\begin{equation}
   \langle w \rangle = \langle \phi_{a_{i_1}} \dots \phi_{a_{i_{|w|}}} \rangle.
\end{equation}
Then Chen's lemma gives the following relation\footnote{In fact, the symbol $\langle \cdot \rangle$ is a character of $(T({\mathscr H}),\shuffle,1_{T({\mathscr H)}})$.}, $\forall \, u , v \in X^*$:
\begin{equation}
\label{Chen}
\left\{
\begin{aligned}  
 \langle u \rangle \langle v \rangle & = \langle u \shuffle v \rangle; \\
 \langle 1 \rangle & = 1
 .
\end{aligned}
\right.
\end{equation}

Some of these iterated integrals have been thoroughly studied, for example the polyzetas: one considers the alphabet $\left\{ x_0 , x_1 \right\}$ and constructs recursively the following integrals: $\forall z \in \C \backslash \left] - \infty , 0 \right] \cup \left[ 1 , + \infty \right[$,
\begin{equation*}
\displaystyle {\rm Li}_{x_0^n}(z) = \frac{\ln^n(z)}{n!},
\end{equation*}

\begin{equation*}
		{\rm Li}_{x_1 w} (z) = \int_0^z \frac{dt}{1-t} {\rm Li}_w(t),
\end{equation*}
and, $\forall w \in X^* x_1 X^*$,
\begin{equation*}
		{\rm Li}_{x_0 w} (z) = \int_0^z \frac{dt}{t} {\rm Li}_w(t). 
\end{equation*}
The specialization of these functions for $z=1$ yields the Multiple Zeta Values (henceforth denoted by MZV) $\zeta({\bf s})$ where the multiindex ${\bf s}$ is obtained from $w$ with the correspondence $w = x_0^{s_1 - 1} x_1 \dots x_0^{s_k - 1} x_1 \leftrightarrow {\bf s} = ( s_1 , \dots s_k )$. One can show that the product of two MZV's is, like the quasi symmetric functions, ruled by the stuffle product $\stuffle$ defined by
\begin{equation}
\begin{array}{rl}
(s_1 , \dots , s_p ) \stuffle ( t_1 , \dots , t_q ) & = \\ 
& s_1 (s_2 , \dots s_p ) \stuffle  ( t_1 , \dots , t_q )\\
+ & t_1 (s_1 , \dots , s_p ) \stuffle ( t_2 , \dots , t_q ) \\
+ & (s_1 + t_1 ) (s_2 , \dots , s_p ) \stuffle ( t_2 , \dots , t_q ).
\end{array}
\end{equation}
Further, coloured polyzetas (\cite{Kreimer,SLC44}) need an indexation by bicompositions $\displaystyle \left( \genfrac{}{}{0pt}{}{s'_1 \dots s'_p}{s_1'' \dots s_p''} \right)$ with a product $\diamond$ given by
\begin{equation}
\begin{array}{rl}
\displaystyle
 \left( \genfrac{}{}{0pt}{}{s'_1 \dots s'_p}{s_1'' \dots s_p''} \right) \diamond \left( \genfrac{}{}{0pt}{}{t'_1 \dots t'_p}{t_1'' \dots t_p''} \right) & = \\
&\displaystyle \left(  \left( \genfrac{}{}{0pt}{}{s'_1}{s_1''} \right) \left( \genfrac{}{}{0pt}{}{s'_2 \dots s'_p}{s_2'' \dots s_p''} \right) \diamond \left( \genfrac{}{}{0pt}{}{t'_1 \dots t'_p}{t_1'' \dots t_p''} \right) \right) \\
+ & \displaystyle \left( \left( \genfrac{}{}{0pt}{}{t'_1}{t_1''} \right) \left( \genfrac{}{}{0pt}{}{s'_1 \dots s'_p}{s_1'' \dots s_p''} \right) \diamond \left( \genfrac{}{}{0pt}{}{t'_2 \dots t'_p}{t_2'' \dots t_p''} \right) \right) \\
+ & \displaystyle \left( \left( \genfrac{}{}{0pt}{}{s_1'+t'_1}{s_1''+t_1''} \right) \left( \genfrac{}{}{0pt}{}{s'_2 \dots s'_p}{s_2'' \dots s_p''} \right) \diamond \left( \genfrac{}{}{0pt}{}{t'_2 \dots t'_p}{t_2'' \dots t_p''} \right) \right) .
 \end{array}
\end{equation}

Even algebras of diagrams $\LDIAG$ (\cite{SLC62}), which need coding with words whose letters, belonging to an alphabet $A$, are composable, are endowed with a product $\uparrow$ of this type. These algebras contain plane bipartite graphs with multiple ordered legs which are in bijection with the elements of $(\mathfrak{MON}^+(X))^*$ where $\mathfrak{MON}^+(X)$ is the set of non void commutative monomials in the variables of the alphabet $X$; formally speaking, 
let $X = \left\{ x_i \right\}_{i \geq 1}$ be an alphabet; denote by $\mathfrak{MON}(X)$ (resp. $\mathfrak{MON}^+(X)$) the monoid of monomials $X^\al$ for $\al \in \N^{(X)}$ (resp. for $\al \in \N^{(X)} \setminus \left\{ 0 \right\}$). Then, the elements of the monoid $(\mathfrak{MON}^+(X))^*$ are \textit{words of monomials} which represent some diagrams.

The bilinear product $\uparrow$ of two diagrams is given on the corresponding words of monomials by
\begin{equation}
\left\{
\begin{aligned}
 1_{(\mathfrak{MON}^+(X))^*} \uparrow w & = w \uparrow 1_{(\mathfrak{MON}^+(X))^*} = w; \\
 a u \uparrow b v & = a (u \uparrow bv) + b (au \uparrow v ) + (a \cdot b) ( u \uparrow v)
 \end{aligned}
\right.
\end{equation}
for all $a, \, b \in \mathfrak{MON}(X)$ and $u, \, v \in (\mathfrak{MON}(X))^*$.

The dualization of the superposition law $(a,b)\rightarrow a \cdot b$ leads to the definition of coproducts given by sums over a semigroup which has the following property: each of its elements has a finite number of decompositions as a product of two elements of the semigroup. This fact motivates the study of such semigroups, called \textit{finite decomposition semigroups}, and of their structure. 

Note that the law of semigroup can be deformed with a bicharacter \cite{Thibon96,Hoffman} or a colour factor \cite{EnjalbertMinh,SLC62,duchamp:hal-00793118}.

The aim of this paper is to present two theorems related to the structure of these semigroups. Section \ref{MotDef} is devoted to the detailed presentation of two examples of the importance of finite decomposition semigroups. In section \ref{DDL}, we give a necessary and sufficient condition for the disjoint direct limit of a family of semigroups to be a finite decomposition semigroup. Finally, in section \ref{finitestruc}, we provide a structure theorem which describes every finite decomposition semigroup as a disjoint direct limit.

\section{Motivations - Definitions}
\label{MotDef}
\noindent 
\subsection{Definitions}
\noindent Let us recall the definition of the finite decomposition property. Let $T$ be a semigroup. We say that $T$ has the \textit{finite decomposition property} (or, equivalently, that $T$ is a finite decomposition semigroup) if, $\forall t \in T$,
\begin{equation}
\label{finitedecomposition}
      \big| \left\{ t_1 , t_2 \in T, \, t_1 \cdot t_2 = t \right\} \big| < \infty.\tag{D}
\end{equation}
We will need the following notation: if $(I,\leq)$ is an ordered set and $\alpha \in I$, then
\[
   \left[ \leftarrow , \alpha \right] = \left\{ \beta \in I , \, \beta \leq \alpha \right\}.
\]
is called the \textit{initial interval} generated by $\alpha$.
\subsection{Motivations}
Our interest for the finite decomposition property comes from the study of several problems of combinatorial physics in which semigroups or monoids with this property are involved\footnote{See \cite{SLC62,duchamp3p} for physics and \cite{duchampint,Hoffman} for the links with number theory.}. We give below two examples.
\subsubsection{Dualizability}
\label{dualizability}
\noindent Another example of the importance of the finite decomposition semigroups comes from the fact that they appear in the study of some bialgebras and, more precisely, in the dualization of the law of the algebra (see, for example \cite{SLC62}; one of the authors recently described some of the features of some semigroup bialgebras \cite{FPSAC2013}). The best-known examples of this kind of bialgebras are given by the shuffle and stuffle algebras where the semigroups involved in the coproduct are respectively the null semigroup and $S = (\N^+ , + )$.

If $k$ is a field and $M$ a semigroup, we denote by $k \left[ M \right]$ the algebra of $M$. It is in duality with itself for the scalar product $\scal{\cdot}{\cdot}$ defined by
\begin{equation*}
 \scal{P}{Q} = \sum_{m \in M} \scal{P}{m} \scal{Q}{m}
\end{equation*}
if $P, \, Q \in k\left[ M \right]$ are polynomials with coefficients $\scal{P}{m}$ and $\scal{Q}{m}$ respectively for all $m \in M$.\\
If $M$ is a finite decomposition semigroup, it is possible to dualize the product of $k \left[ M \right]$. Indeed, one can define the element $\Delta(m) \in k \left[ M \right] \otimes k \left[ M \right]$ by
\[ 
\Delta(m) = \sum_{\genfrac{}{}{0pt}{}{p,q \in M}{pq = m}} p \otimes q
\]
as the sum is finite and then extend $\Delta$ by linearity. One has
\begin{equation*}
\scal{P \cdot Q}{R} = \scal{P \otimes Q}{\Delta(R)}, \, \forall P, Q, R \in k \left[ M \right].
\end{equation*}

\begin{example}
 Let $\mathfrak{MON}^{\rm L}(X) = \left\{ X^\alpha , \, \alpha \in \Z^{(X)} \right\}$ denote the monoid of Laurent monomials. Then the map $\displaystyle \Delta : X^{\alpha} \rightarrow \sum_{\al_1 + \al_2} X^{\al_1} \otimes X^{\al_2}$ takes its values in the large algebra $k \left[ \left[ \mathfrak{MON}^{\rm L}(X) \otimes \mathfrak{MON}^{\rm L}(X) \right] \right]$ since the semigroup of multiindices with values in $\Z$ is not a finite decomposition semigroup.
\end{example}

\subsubsection{Existence of the convolution product}
Let $M$ be a semigroup and ${\mathscr F}(M)$ the space of functions defined on $M$. Assume that $M$ is finite decomposition. Then ${\mathscr F}(M)$ is endowed with the structure of an algebra for the convolution product $\star$ defined by
\[
  (f \star g) (m) = \sum_{m_1 m_2 = m} f(m_1) g(m_2)
\]
for all $f,g \in {\mathscr F}(M)$.

\begin{note}
 In fact, 
 the two motivations above are related. Indeed, the law dual to the coproduct defined in section \ref{dualizability} gives birth to the convolution product of linear forms on $k \left[ M \right]$.
\end{note}

Let $S$ be a finite decomposition semigroup. If $S$ has a neutral $e_S$, we consider $S^{(1)} = S \setminus S^{\times}$ ($S^{\times}$ is the set of invertibles of $S$). The lemmas presented in this paper show that
\begin{itemize}
 \item $S^{1}$ is a sub-semigroup of $S$;
 \item $S^\times$ is a finite group.
\end{itemize}
In the next paragraph, we will see how to iterate this process and reconstruct the initial semigroup from the spare pieces.

\section{Disjoint Direct Limit and Finite decomposition property}
\label{DDL}
\subsection{Disjoint Direct Limit}
Let $(I,\leq)$ be an ordered set. We consider an inductive system of disjoint semigroups $S_\alpha$, $\alpha \in I$. This structure is given by a family of morphisms of semigroups $\phi_{\alpha \beta}: S_\beta \ra S_\alpha$, $\beta \leq \alpha$, which satisfy the following properties:
\begin{subequations} 
\begin{equation} 
\phi_{\alpha \alpha}={\rm Id}_\alpha \text{ for all }\alpha \in I;
\end{equation}
\begin{equation}
\label{comp}
\phi_{\alpha \beta} \circ \phi_{\beta \gamma} = \phi_{\alpha \gamma} \text{ for all }\gamma \leq \beta \leq \alpha \in I.
\end{equation}
\end{subequations}
Then we denote by $S = \underset{\longrightarrow}{\rm DDL}(S_\alpha)$ the \textit{Disjoint Direct Limit} of the system of semigroups which is the semigroup structure on $S = \displaystyle \bigsqcup_{\alpha \in I} S_\alpha$ constructed as follows. Assume that $I$ is a upper half lattice. Then $S$ has the structure of a semigroup for the law $\star$ given by
\begin{equation}
 x \star y = \phi_{\big(\lambda(x) \vee \lambda(y)\big) \lambda(x) } (x) \cdot_{\lambda(x) \vee \lambda(y)} \phi_{\big(\lambda(x) \vee \lambda(y)\big) \lambda(y) } (y)
\end{equation}
where $\lambda(x)$ denotes the unique element of $I$ such that $x \in S_{\lambda(x)}$. Indeed, if $\lambda(x) = \alpha$, $\lambda(y) = \beta$ and $\lambda(z) = \gamma$,
\begin{equation*}
\begin{aligned}
( x \star y ) \star z & = (\phi_{\big(\alpha \vee \beta \big) \alpha } (x) \cdot_{\alpha \vee \beta} \phi_{\big(\alpha \vee \beta \big) \beta } (y) ) \star z \\
& = \phi_{\big( ( \alpha \vee \beta )\vee \gamma \big) (\alpha \vee \beta)} \Big( (\phi_{\big(\alpha \vee \beta \big) \alpha } (x) \cdot_{\alpha \vee \beta} \phi_{\big(\alpha \vee \beta \big) \beta } (y) ) \Big) \cdot_{( \alpha \vee \beta )\vee \gamma} \phi_{\big( (\alpha \vee \beta) \vee \gamma \big) \gamma}(z) \\
& = \Big( \phi_{\big( ( \alpha \vee \beta )\vee \gamma \big) \alpha} (x) \cdot_{(\alpha \vee \beta)\vee \gamma} \phi_{\big( ( \alpha \vee \beta )\vee \beta \big) \alpha} (y) \Big) \cdot_{(\alpha \vee \beta)\vee \gamma} \phi_{\big( (\alpha \vee \beta) \vee \gamma \big) \gamma}(z)
\end{aligned}
\end{equation*}
using the compatibility property \mref{comp} of the morphisms of semigroups $\phi_{\alpha \beta}$. The claim follows from the associativity of the product in $S_{\alpha \vee \beta \vee \gamma}$.

Note that this construction is very similar to the construction of the direct limit of a family of semigroups (which is described, for example, in \cite{B_Set}). It is motivated by the structure of the finite decomposition semigroups (see \ref{finitestruc}): these semigroups can be decomposed as the union of a family of disjoint groups with a finite decomposition semigroup. The disjoint direct limit shows that one can conversely build a semigroup from a family of finite decomposition semigroups.\\
Remark also that in the case where 
\begin{itemize}
   \item the $S_\alpha$'s are monoids with neutral $e_\alpha$;
   \item the morphisms $\phi_{\alpha \beta}$ satisfy $\phi_{\alpha \beta}(e_\beta) = e_\alpha$ for all $\alpha \geq \beta \in I$ (i.e. they are morphisms of monoids);
   \item ${\rm min}(I) = \alpha_0$ exists,
\end{itemize}
then $S$ is a monoid with neutral $e_{\alpha_0}$. Indeed, in that case,
\begin{equation*}
 e_{\alpha_0} \star x = \phi_{\lambda(x) \alpha_0}(e_{\alpha_0}) \cdot_{\lambda(x)} x = e_{\lambda(x)} \cdot_{\lambda(x)} x = x
\end{equation*}
for all $x \in S$.
\subsection{Finite decomposition criterion}
Let $S = {\rm DDL}(S_\alpha)_{\alpha \in I}$ be the disjoint direct limit of a family of semigroups.
\begin{theorem}
The semigroup $S$ is finite decomposition if and only if the following conditions are satisfied:
\begin{enumerate}
 \item[$(i)$] $\forall \alpha \in I$ and $\forall y \in S_\alpha$, $|\left\{ \beta \leq \alpha , \, \phi^{-1}_{\alpha \beta}(y) \neq \emptyset \right\}| < \infty$;
 \item[$(ii)$] every $S_\alpha$ is of finite decomposition type;
 \item[$(iii)$] for all $\alpha \leq \beta \in I$ and for all $x \in S_\alpha$, the fibers $\phi^{-1}_{\alpha \beta}(x)$ of $\phi_{\alpha \beta}$ are finite.
\end{enumerate} 
\end{theorem}

\noindent\textbf{Proof :}
\begin{itemize}
   \item Assume that $S$ is of finite decomposition type. Then it is impossible that one of the intervals $\left\{ \beta \leq \alpha , \, \phi^{-1}_{\alpha \beta}(y) \neq \emptyset \right\}$ be infinite. If that were the case for $y \in S_\alpha$, then $y^2 \in S_\alpha$ has an infinite number of decompositions since $y^2 = y \phi_{\alpha \beta}(x)$ for all $\beta \leq \alpha$ and $x \in \phi^{-1}_{\alpha \beta}(y)$. Moreover, each of the $S_\alpha$'s is a sub-monoid of $S$; hence it is finite decomposition since $S$ is finite decomposition. Finally, the decomposition of $z^2$ presented above also explains why there is no morphism $\phi_{\alpha \beta}$ with an infinite fiber. 
   \item Assume now that the three properties are satisfied. Let $z \in S_{\lambda(z)}$. Consider its decompositions $x y$; they form a set $D_z$ given by $\displaystyle \bigsqcup_{u,v,\alpha,\beta} D(z,u,v,\alpha , \beta)$ where $D(z,u,v,\alpha,\beta) = \left\{ (x,y) \in S_\alpha \times S_\beta \text{ such that } \phi_{\lambda(z)\alpha}(x)=u \text{ and } \phi_{\lambda(z)\beta}(y)=v \right\}$. Remark that in order that $D(z,u,v,\alpha,\beta) \neq \emptyset$, one must have $\alpha \vee \beta = \lambda(z)$. There is a finite number of $\alpha$ and $\beta$ such that $\alpha < \lambda(z)$ and $\beta < \lambda(z)$ because of the structure of $I$. Moreover, there is also a finite number of $x$ and $y$ respectively in $S_\alpha$ and $S_\beta$ such that $u = \phi_{\lambda(z) \alpha}(x)$ and $v = \phi_{\lambda(z) \beta}(y)$ because the morphisms are finite fibers. Finally, since the $S_\alpha$'s are finite decomposition, there is a finite number of decompositions of $z$ as a product of elements of $S_{\lambda(z)}$. 
Hence $S$ is finite decomposition.\cqfd
\end{itemize} 

\begin{example}
The following shows an example of ${\rm DDL}$ which is of finite decomposition type but whose intervals $\left[ \leftarrow , \alpha \right]$ are all infinite. It proves that the first condition can not be replaced by the finiteness of every interval $\left[ \leftarrow , \alpha \right]$.\\
Let $S_k = \left[ k , +\infty \right[$ be the additive semigroup of integers greater than or equal to $k$. The disjoint direct limit of the family $(S_k)_{k \geq 0}$ is the semigroup $S = \left\{ (k,x) \in \N^2, \, k \leq x \right\}$ with product $\star = \wedge_\N \times +$ where
\begin{equation*}
 (k_1 , y_1) \star (k_2 , y_2 ) = ( k_1 \wedge_\N k_2 , y_1 + y_2). 
\end{equation*}
The morphisms $S_k \stackrel{\phi_{k\ell}}{\longrightarrow} S_\ell$ associate $(\ell,y)$ to $(k,y)$ for $\ell \leq k$ (see Fig. \ref{Fig}). The set of index of the semigroups is $\N$ with an order $\prec$ such that
\[ 
   \alpha \prec \beta \Leftrightarrow \alpha \wedge \beta = \beta.
\]
Hence, the intervals $\left[ \leftarrow , \alpha \right] = \left\{ \beta , \beta \wedge \alpha = \alpha \right\} = \left[ \alpha , + \infty \right[$ are all infinite.
\begin{center}
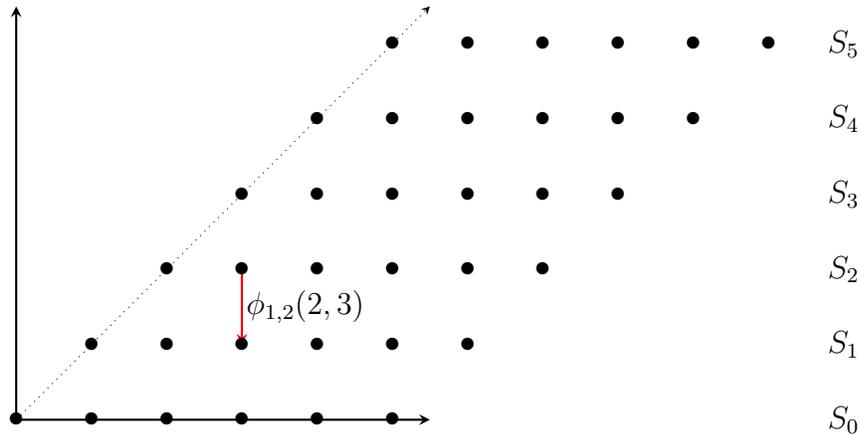
\begin{figure}[h]
\begin{equation*}
\begin{tikzpicture}[>=stealth]
\draw[->,thick] (0,0) -- (5.5,0);
\draw[->,thick] (0,0) -- (0,5.5);
\draw[->,dotted] (0,0) -- (5.5,5.5);

\draw[->, thick, color=red] (3,2) -- (3,1);

\draw (3.7,1.5) node {$\,\,\,\,\phi_{1,2}(2,3)$};

\draw (11,0) node {$S_0$};
\draw (11,1) node {$S_1$};
\draw (11,2) node {$S_2$};
\draw (11,3) node {$S_3$};
\draw (11,4) node {$S_4$};
\draw (11,5) node {$S_5$};

\draw (0,0) node {$\bullet$};
\draw (1,0) node {$\bullet$};
\draw (2,0) node {$\bullet$};
\draw (3,0) node {$\bullet$};
\draw (4,0) node {$\bullet$};
\draw (5,0) node {$\bullet$};

\draw (1,1) node {$\bullet$};
\draw (2,1) node {$\bullet$};
\draw (3,1) node {$\bullet$};
\draw (4,1) node {$\bullet$};
\draw (5,1) node {$\bullet$};
\draw (6,1) node {$\bullet$};

\draw (2,2) node {$\bullet$};
\draw (3,2) node {$\bullet$};
\draw (4,2) node {$\bullet$};
\draw (5,2) node {$\bullet$};
\draw (6,2) node {$\bullet$};
\draw (7,2) node {$\bullet$};

\draw (3,3) node {$\bullet$};
\draw (4,3) node {$\bullet$};
\draw (5,3) node {$\bullet$};
\draw (6,3) node {$\bullet$};
\draw (7,3) node {$\bullet$};
\draw (8,3) node {$\bullet$};

\draw (4,4) node {$\bullet$};
\draw (5,4) node {$\bullet$};
\draw (6,4) node {$\bullet$};
\draw (7,4) node {$\bullet$};
\draw (8,4) node {$\bullet$};
\draw (9,4) node {$\bullet$};

\draw (5,5) node {$\bullet$};
\draw (6,5) node {$\bullet$};
\draw (7,5) node {$\bullet$}; 
\draw (8,5) node {$\bullet$};
\draw (9,5) node {$\bullet$};
\draw (10,5) node {$\bullet$};



 \end{tikzpicture}
\end{equation*}
\caption{Illustration of the ${\rm DDL}$ of the family $(S_k)_{k \geq 0}$.}
\label{Fig}
\end{figure}
\end{center}
\end{example}

\section{Structure of the finite decomposition monoids}
\label{finitestruc}
\subsection{Structure theorem}
\label{secth}
 Let $T$ be a semigroup. One defines two sequences $(T_n)_{n \in \N}$ and $(D_n)_{n \in \N}$ by
 \begin{itemize}
  \item $T_0=T, \qquad D_0=\emptyset$;
  \item if $T_n,\, D_n$ are constructed, 
  \begin{itemize}
  \item either $T_n$ has no neutral and we stop with $D=D_{n}$;
  \item or $T_n$ has a neutral $e_n$ and 
\begin{equation}
	T_{n+1}=T_n\setminus T_n^{\times}\ ;\ D_{n+1}=D_n\cup \{n+1\}.
\end{equation}
 \end{itemize}
\end{itemize}
For convenience, we denote by $G_n$ the group of invertible elements of $T_n$ : $G_n=T_n^{\times}$ whenever $T_n$ admits a neutral $e_n$.

\begin{theorem}
 Let $T$ be a finite decomposition semigroup. Then 
\begin{enumerate}
\item[$i)$] either $D = \left\{ 1 , \dots , N \right\}$ is finite and $T=(\displaystyle \bigsqcup_{1\leq n\leq N}G_n)\sqcup T_{N+1}$ where $T_{N+1}$ is a finite decomposition semigroup without neutral; $T_m=(\displaystyle \bigsqcup_{m\leq n\leq N}G_n)\sqcup T_{N+1}$ is a sub-semigroup of $T$ and there exists a family of morphisms of monoids $\phi_{ij}~: G_j\ra G_{i};\ \phi_{Ni}~: G_i\ra  T_{N}$.\\
\item[$ii)$] or $D$ is infinite and $T=(\displaystyle \bigsqcup_{n \geq 0}G_n)$; $T_m=(\displaystyle\bigsqcup_{m\leq n\leq N}G_n)$ is a sub-semigroup of $T$ and there exists a family of morphisms of monoids $\phi_{ij} : G_j\ra G_{i}$.\\
\end{enumerate}
\end{theorem}

\subsection{Lemma 1}
We will need the following lemma for the proof of the theorem.
\begin{lemma}
\label{lemme1}
 Let $T$ be a finite decomposition semigroup with unit $1$. Then the following properties are equivalent:
 \begin{enumerate}
  \item[$(i)$] $u \in T$ is right invertible;
  \item[$(ii)$] $u$ is left invertible;
  \item[$(iii)$] $u$ is cyclic (and hence invertible).
 \end{enumerate}
\end{lemma}
\noindent\textbf{Proof: }Let $u \in T$ be a right invertible element. Then, for all $n \in \N$, $u^n v^n = 1$. Since $T$ has the finite decomposition property, it is impossible that all the decomposition of $1$ of the form $u^n v^n$ are different. Thus, there exists $p>0$ such that $(u^{n+p},v^{n+p}) = (u^n , v^n)$. Hence $u^{n+p} = u^n$ and one has
\begin{equation}
 1 = u^n v^n = u^{n+p} v^n = u^p
\end{equation}
which proves that $u$ is invertible since $u u^{p-1} = 1 = u^{p-1} u$.

The same proof holds if $u$ is left invertible and both $(iii) \Rightarrow (i)$ and $(iii) \Rightarrow (ii)$ are trivial, hence the claim.\cqfd

\subsection{Lemma 2}
In this section, we use the notation of section \ref{secth}. Let $T$ be a finite decomposition semigroup.
\begin{lemma}
\label{lemme2}
For all $x\in T$ and $n\in D$, one has
\begin{equation}
	e_nxe_n=e_nx=xe_n \in T_n.
\end{equation}
\end{lemma}
\noindent\textbf{Proof: }Let $x \in T$ and $n \in D$. We denote by $i_0 = \underset{n}{{\rm max}} \left\{ xe_n \in T_n \right\}$. If $i_0 = N$ we are done. Assume that $i_0 < N$. One has $x e_{i_0 +1} = x e_{i_0} e_{i_0 + 1}$ (since $e_{i_0}$ is the neutral of $T_{i_0}$ which contains $e_{i_0+1}$); $x e_{i_0} \in T_{i_0}$ hence $x e_{i_0 +1} \in T_{i_0}$. But $x e_{i_0 +1} \notin G_{i_0}$ (if that were the case, then lemma \ref{lemme1} would imply that $e_{i_0+1}$ belongs to $G_{i_0}$; this is impossible by definition of $e_{i_0+1}$); thus $x e_{i_0+1} \in T_{i_0} \setminus G_{i_0} = T_{i_0+1}$. This is not possible by definition of $i_0$. Necessarily $i_0 = N$. The same argument proves that $j_0 = \underset{n}{{\rm max}} \left\{ e_n x \in T_n \right\}=N$. The claims follows from the fact that $e_n$ is the neutral of $T_n$; hence $x e_n = e_n (x e_n) = (e_n x ) e_n = e_n x$. \cqfd

\subsection{Proof of the theorem}

Note that, for all $n \in D$, $T_n$ is a semigroup: if $x$ and $y$ belong to $T_n$, $xy$ belong to $T_n$. Indeed, if that is not the case, $xy$ belongs to $G_{n-1}$; then $x$ and $y$ are invertible but it is not possible since $x$ and $y$ belong to $T_n=T_{n-1} \setminus G_{n-1}$. 

Lemma \ref{lemme2} implies that for all $x \in T$ and for all $n \in D$, $e_n x e_n \in T_n$. Hence, $\phi_n : T \ra T_n$ defined by $\phi(x) = e_n x e_n$ is a morphism of monoids (since $\phi_n(xy) = e_n x y e_n = e_n x e_n e_n y e_n = \phi_n(x) \phi_n(y)$; of course $\phi_n(e_j) = e_n$). The restrictions $\phi_{i}\Big|_{G_j} : G_j \ra G_i$ define the morphisms $\phi_{ij}$.

From now on, we assume that $T$ is a finite decomposition semigroup. As an intersection of a non empty (as soon as $e_1$ exists) family of semigroups, $T_{N+1} = T \displaystyle \setminus \bigsqcup_{n \in D} G_n$ is a finite decomposition semigroup. \\ 
Assume that $D$ be infinite; then $T \displaystyle \setminus \left( \bigsqcup_{n \in D} G_n \right)$ is empty. Indeed, if there were $t \in T \displaystyle \setminus \left( \bigsqcup_{n \in D} G_n \right)$, then $t \in T_n$ for all $n \in D$; thus $e_n t = t$ for all $n \in D$ and $t$ is an element of $T$ that has an infinite number of decompositions; this is not possible. \\
If $D$ is finite, $T_{N+1}$ has no neutral. Indeed, assume that there be a neutral $e_{N+1} \in T_{N+1}$. Then $e_{N} e_{N+1} = e_{N}$; hence $e_{N+1}$ is left invertible and thus invertible in $T_N$; this is not possible since $e_{N+1}$ belongs to $T_{N+1} = T_{N} \setminus G_N$.

\section{Conclusion}
We have illustrated the importance of finite decomposition semigroups for the computation of different generalized stuffle products. \\
It turns out that every finite decomposition semigroup is the disjoint direct limit of finite groups and possibly a finite decomposition semigroup without neutral.\\
Moreover, it is possible to apply the disjoint direct limit process to construct new finite decomposition semigroups.

\section*{Acknowledgements}
The authors take advantage of these lines to acknowledge support from the French Ministry of Science and Higher Education under Grant ANR PhysComb and local support from the Project "Polyzetas". They also wish to express their gratitude to K. Penson and C. Tollu for fruitful discussions.

\bibliographystyle{alpha}
\bibliography{Biblio}
\end{document}